\documentclass[11pt]{article}

\usepackage{amscd, mathtools, amsmath, amssymb, fancyhdr, color}

\usepackage{fourier}
\usepackage{libertine}
\usepackage[backref=page]{hyperref}
\renewcommand*{\backref}[1]{}
\renewcommand*{\backrefalt}[4]{%
    \ifcase #1 (Not cited.)%
    \or        (Cited on page~#2.)%
    \else      (Cited on pages~#2.)%
    \fi}

\newcommand{\version}{version 2.0,\ \  06.01.2017}

\numberwithin{equation}{section}

\newcommand{\ZZ}{\mathbb{Z}}

\newcommand{\la}{\lambda}

\def\eqref#1{(\ref{#1})}

\newcommand{\arrow}{{\:\longrightarrow\:}}
\newcommand{\Z}{{\Bbb Z}}
\newcommand{\C}{{\Bbb C}}
\newcommand{\R}{{\Bbb R}}
\newcommand{\Q}{{\Bbb Q}}
\renewcommand{\H}{{\Bbb H}}

\def\1{\sqrt{-1}\:}
\newcommand{\restrict}[1]{{\left|_{{\phantom{|}\!\!}_{#1}}\right.}}
\newcommand{\cntrct}                
{\hspace{2pt}\raisebox{1pt}{\text{$\lrcorner$}}\hspace{2pt}}

\newcommand{\calo }{{\cal O}}

\renewcommand{\tilde}{\widetilde}
\renewcommand{\bar}{\overline}
\renewcommand{\phi}{\varphi}
\renewcommand{\epsilon}{\varepsilon}
\renewcommand{\geq}{\geqslant}

\newcommand{\Hom}{\operatorname{Hom}}

\newcommand{\rk}{\operatorname{rk}}

\newcommand{\Stab}{\operatorname{Stab}}
\newcommand{\reg}{\operatorname{reg}}

\newcommand{\ie}{{\em i.e. }}
\newcommand{\eg}{{\em e.g. }}

\renewcommand{\Im}{\operatorname{Im}}
\renewcommand{\rk}{\mathrm{rank}}

\def\ra{\arrow}

\newcounter{Mycounter}[section]
\newcounter{lemma}[section]
\newcounter{claim}[section]
\renewcommand{\theclaim}{{Claim \thesection.\arabic{claim}}}
\newcommand{\claim}{%
     \setcounter{claim}{\value{Mycounter}}
     \refstepcounter{claim}
     \stepcounter{Mycounter}
     { \bf \theclaim.\ }}

\newcounter{sublemma}[section]
\newcounter{corollary}[section]
\newcounter{theorem}[section]
\renewcommand{\thetheorem}{{Theorem \thesection.\arabic{theorem}}}
\newcommand{\theorem}{%
     \setcounter{theorem}{\value{Mycounter}}
     \refstepcounter{theorem}
     \stepcounter{Mycounter}
     { \bf \thetheorem.\ }}
     
     \newcounter{theoremm}[section]
\newcounter{conjecture}[section]
\newcounter{proposition}[section]
\renewcommand{\theproposition}
       {{Proposition \thesection.\arabic{proposition}}}
\newcommand{\proposition}{%
     \setcounter{proposition}{\value{Mycounter}}
     \refstepcounter{proposition}
     \stepcounter{Mycounter}
     { \bf \theproposition.\ }}

\newcounter{definition}[section]
\renewcommand{\thedefinition}
       {{Definition~\thesection.\arabic{definition}}}
\newcommand{\definition}{%
     \setcounter{definition}{\value{Mycounter}}
     \refstepcounter{definition}
     \stepcounter{Mycounter}
     { \bf \thedefinition.\ }}

\newcounter{example}[section]
\renewcommand{\theexample}{{Example \thesection.\arabic{example}}}
\newcommand{\example}{%
     \setcounter{example}{\value{Mycounter}}
     \refstepcounter{example}
     \stepcounter{Mycounter}
     { \bf \theexample:\ }}

\newcounter{remark}[section]
\renewcommand{\theremark}{{Remark \thesection.\arabic{remark}}}
\newcommand{\remark}{%
     \setcounter{remark}{\value{Mycounter}}
     \refstepcounter{remark}
     \stepcounter{Mycounter}
     { \bf \theremark.\ }}

\newcounter{problem}[section]
\newcounter{question}[section]
     \newcounter{fact}[section]
\renewcommand{\leftmark}%
{{\scriptsize  Flat affine subvarieties in OT manifolds}}

\makeatletter

\@addtoreset{equation}{section}
\@addtoreset{footnote}{section}
\makeatother

\def\blacksquare{\hbox{\vrule width 5pt height 5pt depth 0pt}}
\def\endproof{\blacksquare}

\begin{document}
\begin{center}
{\LARGE\bf Flat affine subvarieties in Oeljeklaus-Toma manifolds
}\\[.1in]

Liviu Ornea\footnote{Partially supported by a grant of Ministry of Research and Innovation, CNCS - UEFISCDI,
project number PN-III-P4-ID-PCE-2016-0065, within PNCDI
III.}
, Misha Verbitsky\footnote{Partially supported by the 
Russian Academic Excellence Project '5-100'.\\[1mm]

}, and Victor Vuletescu${}^1$
\\[1mm]

\end{center}

{\small
\hspace{0.15\linewidth}
\begin{minipage}[t]{0.7\linewidth}
{\bf Abstract:} \\ 
The Oeljeklaus-Toma (OT-) manifolds are compact, 
complex, non-K\"ahler manifolds
constructed by Oeljeklaus and Toma, and generalizing the Inoue surfaces.
Their construction uses the number-theoretic data:
a number field $K$ and a torsion-free subgroup $U$
in the group of units of the ring of integers of $K$, with
rank of $U$ equal to the number of real embeddings of $K$.
OT-manifolds are equipped with a torsion-free flat affine connection
preserving the complex structure (this structure is known as
``flat affine structure'').
We prove that any complex subvariety 
of smallest possible positive dimension 
in an OT-manifold is also flat affine.
This is used to show that if all 
elements in $U\setminus\{1\}$ are 
primitive in $K$, then $X$ contains no proper 
analytic subvarieties.
\end{minipage}
}

\tableofcontents

\section{Introduction}

The OT (Oeljeklaus-Toma) manifolds were discovered by
K. Oeljeklaus and M. Toma in 2005 (\cite{ot}). 
These manifolds are complex solvmanifolds 
generalizing the Inoue surfaces of class S${}^0$
(\cite{inoue}). 

The construction of 
Oeljeklaus and Toma is based on number-theoretic date.
However, the geometry of OT-manifold is best understood
using the Lie group theory. 

Let $G$ be a Lie group equipped with a right-invariant
integrable complex structure.
Recall that a {\bf group manifold}
is the quotient $G/\Gamma$ of $G$ by the right
action of a discrete, cocompact subgroup $\Gamma\subset G$.
A {\bf complex solvmanifold} is a group manifold with a solvable
group $G$. The notion of a solvmanifold is due to 
G. D. Mostow, who proved a structure theorem
for (real) solvmanifolds in his first paper
\cite{_Mostow:solvmanifolds_}. 
The corresponding notion of a {\it complex} solvmanifold
is probably due to K. Hasegawa 
\cite{_Hasegawa_} who also classified
2-dimensional complex solvmanifolds.

In the case of Oeljeklaus-Toma manifolds,
the solvable Lie group is obtained as follows.
Recall that a {\bf metabelian} group is a semidirect
product of two abelian groups. Consider
two abelian Lie groups $A_\R$ and $H_\R$
associated with a number field $K$.
We define $A_\R:= \calo_K \otimes_Z \R$ and
$H_\R:= U \otimes_\Z \R$, where
$U$ is a free abelian subgroup
in the group $\calo_K^*$ of units in 
the integers ring $\calo_K$ of $K$.
There is a natural action of
$U$ on $\calo_K$, allowing one
to define the semidirect product
$G:=A_\R \ltimes H_\R$. The corresponding
cocompact discrete group
is $\Gamma:=\calo_K^+\ltimes U$, where $\calo_K^+$
is the additive group of $\calo_K$.
The OT-manifold is $G/\Gamma$,
with the right-invariant complex 
structure defined explicitly in Section \ref{secOT}.

OT-manifolds provide  a counterexample
to Vaisman's conjecture \cite{vai} 
which was open for 25 years.
They are non-K\"ahler flat affine
complex manifolds\footnote{A complex manifold is
called {\bf flat affine} if it is equipped with a
flat torsion-free connection preserving the
complex structure.} of algebraic dimension 0
(\cite{ot}). Since their discovery in
2005, OT-manifolds were the subject of
much research of complex geometric
and number theoretic nature (\cite{bo}, \cite{ov_mrl_11},
\cite{mt}, \cite{sima}, \cite{_Parton_Vuletescu_}, \cite{_Braunling_}).

It is known that OT-manifolds have no complex curves,
\cite{sima}, and for $t=1$ they have no complex
subvarieties (see \cite{ov_mrl_11}, where the proof makes
explicit use of the LCK structure). Moreover, 
all surfaces contained in OT-manifolds are
blow-ups of Inoue surfaces $S^0$,
\cite{_sima:surfaces_}. However, in general, there is no
characterization of the possible subvarieties of OT
manifolds. The aim of this paper is to give a sufficient
condition for an OT-manifold to not have
submanifolds. In  Section \ref{secOT}, we prove:

\hfill

\theorem\label{mainn} Let $X=X(K, U)$ be an 
OT-manifold. Assume that any element $u\in U\setminus\{1\}$
is a primitive element for the number field 
$K.$ Then $X$ contains no proper 
complex analytic subvarieties.

\hfill

We also prove the following theorem.
Recall that {\bf a flat affine manifold}
is a manifold equipped with a torsion-free flat 
connection. By construction,  OT-manifolds come
equipped with a flat affine structure. 
A submanifold of $Z\subset M$ of a flat affine manifold
is called {\bf flat affine}
if locally around any smooth point $z\subset Z$,
the sub-bundle $TZ\subset TM\restrict Z$ is preserved
by the flat affine connection. Notice that all
flat affine manifolds are equipped with 
local coordinates such that the transition
functions are affine, and in these coordinates
$Z$ is an affine subspace.

\hfill

\theorem\label{_flat_subva_Theorem_}
Let $M$ be an OT-manifold, and $Z\subset M$
an irreducible complex subvariety of smallest possible
positive dimension. Then $Z$ is a smooth flat affine
submanifold of $M$.

{\bf Proof:} See \ref{_subvarieties_flat_Remark_}. \endproof

\hfill

The paper is organized as follows. Section 2 will describe
the construction and main properties of OT-manifolds,
Section 3 provides examples of OT submanifolds, in Section
4 we prove that all holomorphic maps from tori to OT
manifolds are constant, while in Section 5 we give the
proof  of \ref{mainn}.

\section{OT-manifolds}

We briefly describe  the construction of Oeljeklaus-Toma manifolds, following \cite{ot}. 

Let $K$ be a number field which has $2t$
complex embeddings denoted $\tau_i, \bar \tau_i$ and 
$s$ real ones denoted $\sigma_i$,  $s>0$, $t>0$ (for what needed in this paper about number theory, see {\em e.g.} \cite{milne2}).  

Denote $\calo_K^{*,+}:=\calo_K^*\cap \bigcap_i
\sigma^{-1}_i(\R^{>0})$.
Clearly, $\calo_K^{*,+}$ is a finite index
subgroup of the group of units of $\calo_K$.

 Let $\H=\{y\in \C\,;\, \Im y>0\}$ be the upper half-plane. 
For any $\zeta\in \calo_K$  define the
automorphism $T_\zeta$  of $\H^s\times \C^t$  by
$$T_\zeta(x_1,\ldots,x_t, y_1,\ldots,
y_{s})=\bigg(x_1+\tau_1(\zeta),\ldots,
x_t+\tau_t(\zeta),y_1+\sigma_1(\zeta),\ldots,
y_{s}+\sigma_{s}(\zeta)\bigg).
$$
Similarly, for any totally positive\footnote{An element of
  a number field is
  called ``totally positive'' if it is mapped to a
  positive number under all real embeddings.} 
unit $\xi$, let $R_\xi$ be the automorphism of $\C^t\times \H^s$ defined by
$$R_\xi(x_1,\ldots,x_t, y_1,\ldots, y_{s})=\bigg (\tau_1(\xi)x_1,..., \tau_t(\xi)x_t,\sigma_1(\xi) y_1,  ..., \sigma_t(\xi) y_s\bigg).$$
Note that the totally positivity of $\xi$ is needed for $R_\xi$ to act on $\H^s\times \C^t.$

For any subgroup $U\subset \calo_K^{*, +}$, the above maps define a free action of the
semidirect product $U\ltimes \calo_K$ on $\H^s\times
\C^t$.

It is proven in \cite{ot} that one can always find {\bf admissible subgroups} $U$ such that the above action is discrete and cocompact.  Note that if $U$ is an admissible subgroup then necessarily one has 
$$\rk_\ZZ(U)+\rk_\ZZ(\calo_K)=2(s+t),$$
and hence $\rk_\ZZ(U)=s.$ This explains why the condition $t>0$ is needed: otherwise we would have  $\rk\ \calo_K^*=s-1,$ strictly less than $s,$ and thus admissible subgroups could not exist.

\hfill

\definition { (\cite{ot})}\label{defot} For an  admissible subgroup $U$, the quotient
$X(K,U):=({\H}^s\times \C^t)/(U\ltimes \calo_K)$ is called 
 an {\bf Oeljeklaus-Toma manifold} (OT-manifold for short). 

\hfill

\remark\label{covering} It was observed in \cite{mt} that
in the previous construction one may take instead of the
ring of integers $\calo(K)$ {\em any} (additive) subgroup
$H\subset (K, +)$ which equals $\calo (K)$ up to finite
index, \ie either $H\subset \calo (K)$ or $\calo
(K)\subset H$ with finite index. We let $X(K,H,U)$ the
resulting manifold. Note that the OT-manifolds in
\ref{defot} correspond to the case $H=\calo_K$. Clearly,
any such $X(K,H,U)$ is isogeneous to $X(K,U)$, that is,
$X(K,H,U)$ is a finite cover of a finite quotient of $X(K,U)$.

\hfill

\remark\label{_inoue_dim_1_Remark_}
For $s=t=1$, one recovers 
a version of the classical construction
used by Inoue to define the 
Inoue surfaces of class $S^0$ (\cite{inoue}). 
In \cite{inoue}, no number theory was employed. However,
the matrix $M\in \mathrm{SL}(3, \Z)$ used in \cite[\SS 2]{inoue} to construct the
Inoue surface of class $S^0$ gives a cubic number field,
generated by its root, and this field can be used to 
recover $M$ in a usual way. If one applies the
Oeljeklaus-Toma construction to this cubic field,
one would obtain the Inoue surface associated with $M$.

\hfill

\remark\label{_OT_LCK_Remark_}
All OT-manifolds ($s>0$) are non-K\"ahler, but for $t=1$
they admit locally conformally K\"ahler (LCK) metrics (see
\cite{do} for this notion).

\hfill

\definition\label{simple} (\cite{ot}) An OT-manifold is
called {\bf of simple type},  if $U\not\subset\Z$ and it
does not satisfy any of 
the following equivalent\footnote{Equivalence follows from
\cite[Lemma 1.4]{ot}.} conditions:
\begin{enumerate}
\item The action of $U$ on $\calo_K$  admits a proper, non-trivial, invariant submodule of lower rank.
\item There exists some proper, intermediate field extension $\Q\subset K'\subset K$, with $U\subset\calo^*_{K'}$.
\end{enumerate}

\remark\label{subv} a) A simple type OT-manifold has no proper OT submanifolds with the same group of units $U$. Also, note the difference towards the notion of simplicity in \cite{cdv}.

b) {If $K$ is a number field, then, by
  Dirichlet's theorem, under the logarithmic embedding its
  group of units identifies (up to its subgroup of roots
  of unity) with a full lattice in a real vector space. Denote
  this space by $V_K.$ 
Similarly, if $K'\subset K$ is some subfield of $K$, then the group of units of $K'$ identifies  with a lattice in a proper vector subspace $V_{K'}\subset V_K.$ As $K$ has finitely many subfields, and as the admissible  group of units $U$ of an OT-manifold can be chosen generically in the group of units of $K$, we see that the  OT-manifolds of simple type are generic.}

\hfill

\section{Examples of submanifolds in  OT-manifolds}

A simple example of an OT-manifold   embedded in a larger OT-manifold which is {\em not of simple type} in the sense of Oeljeklaus-Toma is constructed in \cite[Remark 1.7]{ot}.  

We now provide an example of an OT submanifold embedded in an OT-manifold which is {\em  of simple type} in the sense of Oeljeklaus-Toma.

\hfill

\example\label{expl2} Take $L=\Q[X]/(X^3-2)$; then $L$ has one real embedding $\tau_1$ and $2t=2$ complex ones $\tau_2, \tau_3(=\overline{\tau}_2).$  Note that $U_L=\calo_L^{*,+}$  is a free group of rank one, and denote $u_1$ be a generator for $U_L$. Then $U_L$  is an admissible group, and  let $S=X(L, U_L)$ is the corresponding OT-manifold (an Inoue surface $S^0$).

Now take $K=\Q[X]/(X^6-2).$ The field $K$ is an extension of degree $2$ of $L$ which has two real embeddings $\sigma_1, \sigma_2$ (which both extend the embedding $\tau_1$ of $L$) and four complex embeddings: $\sigma_3, \sigma_4$ (which extend $\tau_2$) and $\sigma_5=\overline{\sigma}_3, \sigma_6=\overline{\sigma}_4$ (which extend $\tau_3=\overline{\tau}_2).$
Consider the unit $u_2\in \calo_K^{*+}$ such that $\sigma_1(u_2)=(\sqrt[6]{2}-1)^2.$ Then $\sigma_2(u_2)=(\sqrt[6]{2}+1)^2$, and  hence  the subgroup $U_K\subset \calo_K^{*+}$ generated by $u_1$ and $u_2$  is admissible, since the projection on the first two factors of their logarithmic embedding is
$$
\left(
\begin{array}{ccc}
\log(u_1) & \log(u_1)\\[.1in]
2\log(\sqrt[6]{2}-1) & 2\log(\sqrt[6]{2}+1)
\end{array}
\right)
$$
which is of maximal rank.

Let $X=X(K, U_K)$ be the corresponding OT-manifold, $X=\H^2\times \C^2/(U_K\ltimes \calo_K).$ Define the map 
$i:S\ra X$ by $$i([w, z])=[w,w,z,z],$$
where we denoted by $[x]$ the equivalence classe of $x$.
Clearly, $i$ is well-defined.

\hfill

\claim The map $i$ is injective. 

Indeed, if $i([w, z])=i([w', z']$ then there exists $u\in U_K, a\in \calo_K$ such that
\begin{equation*}
\begin{split}
w&=\sigma_1(u)w'+\sigma_1(a),\\
w&=\sigma_2(u)w'+\sigma_2(a).\\
\end{split}
\end{equation*}
This implies 
$$(\sigma_1(u)-\sigma_2(u))w'=\sigma_2(a)-\sigma_1(a).$$  
If $\sigma_1(u)\not=\sigma_2(u)$, then $w'\in \R$, which is not possible, and hence $\sigma_1(u)=\sigma_2(u)$. This yields $u\in L$, thus $u\in U_L$. Moreover, $\sigma_1(a)=\sigma_2(a)$, and hence  also $a\in \calo_L$. But then  $[w, z]=[w', z']$. \endproof

\hfill

\remark The above constructed $X$ is of simple type.  Indeed, the unit $u_2$ is a primitive element for $K$, hence there is no proper subfield $K'\subset K$ containing $u_2.$
\endproof



\section{Holomorphic maps from and to tori}

In the proof of the main result we shall need the following result interesting in itself:

\hfill

\proposition\label{tori}
Let $X$ be an  OT-manifold and $T$ a complex torus. Then any holomorphic map $f:T\ra X$ must be constant.

\hfill

{\bf Proof:} Let $T=\C^d/\Lambda$, and  $X=\H^s\times \C^t/(U\ltimes A)$. Note that  $\pi_1(T)=\Lambda$,  and $\pi_1(X)=U\ltimes A$. Let $I$ be the image of the natural morphism $f_*:\pi_1(T)\ra \pi_1(X).$ With the above identifications, we let:
 $$f_*(\lambda)=\gamma_\lambda=(u_\lambda, a_\lambda),\quad \la\in\Lambda.$$ 
 Let now   $\tilde{f}:\tilde{T}\ra \tilde{X}$ be a lift of $f$ at the universal covers. We then have:
\begin{eqnarray}\label{univ} 
\tilde{f}(t+\lambda)=\gamma_\lambda(\tilde{f}(t)), \forall t\in \tilde{T}.
\end{eqnarray}
Let $f_1=\mathrm{pr}_1\circ \tilde{f}$, where $\mathrm{pr}_1:\H^s\times \C^t$ is the projection onto the first factor. Then $f_1$ is a map from $\C^d$ into $\H$, and hence by Liouville's theorem  it must be constant. It follows that   the first component of the map $\tilde{f}$ is a constant, say $w_1$. Then \eqref{univ} implies $$w_1=\sigma_1(u_\lambda)w_1+\sigma_1(a_\lambda).$$ 
Now if $\sigma_1(u_\lambda)\not=1$ for some $\lambda\in \Lambda$, then $w_1$ would be real, which is impossible.  It follows that $u_\lambda=1$ for all $\lambda\in \Lambda$ and thus $I$ is actually a subgroup of $A$.

But then we have
$w_1=w_1+\sigma_1(a_\lambda)$, and hence $\sigma_1(a)=0$, yielding $a_\lambda=0$ for all $\lambda\in \Lambda.$
 This implies $I=\{0\}$, and thus $\tilde{f}$ factors through a map from $T$ to the universal cover $\tilde{X}$ of $X,$ which is constant as $T$ is compact.
\endproof

\hfill

\remark It is easy to see that  conversely, every holomorphic map from an OT-manifold to a torus must be constant. This is because the Albanese torus of an OT-manifold is trivial, since there are no non-zero closed holomorphic 1-forms on an OT-manifold (\cite[Proposition 2.5]{ot}).

\section{The proof of the  main result}\label{secOT}

\theorem\label{main} Let $X=X(K, U)$ be an OT-manifold. Assume that any element $u\in U\setminus\{1\}$ is a primitive element for $K.$ Then $X$ contains no proper analytic subspaces.

\hfill

{\bf Proof.}  We argue by contradiction. Let $Z\subset X$ be an analytic connected proper subpspace of minimum positive dimension. By a result of S.M. Verbitskaya, an OT-manifold cannot contain curves, and hence  $\dim(Z)\geq 2$ (note that the proof in \cite{sima} can be easily  extended to cover the singular case, too). Moreover, as $Z$ is of minimum positive dimension, we deduce that $Z$ contains no Weil divisors, and it has at most finitely many isolated singularities.

Let $Z_{\reg}=Z\setminus \mathrm{Sing}(Z)$ be the regular
part of $Z$. Then the Remmert-Stein's theorem (see \eg 
\cite[Theorem 6.9, p. 150]{fg})  implies that $Z_{\reg}$ has
no divisors.

For any $i=1,\ldots, \dim(X)$ let $L_i$ be the flat line bundle on $X$ associated with the representation $\varrho_i:\pi_1(X)\ra \C^*$, $\varrho_i(u, a)=\sigma_i(u)$ (here we identified $\pi_1(X)$ with $U\ltimes \calo_K$). Then $L_i$ is locally generated by $\frac{\partial}{\partial z_i}$ and the tangent bundle $T_X$ is naturally identified with the direct sum
$$T_X=\bigoplus_{i=1}^n L_i.$$

We want to understand the restriction of $T_X$ to $Z$. It will be enough to look at the regular part $Z_{\reg}$ on which  we have the exact sequence:
\begin{equation}\label{seq}
0\ra T_{Z_{\reg}}\stackrel{i}{\ra} \bigoplus_{i=1}^n {L_i}\restrict{ Z_{\reg}}\ra {\cal N}_{Z_{\reg}\vert X}\ra 0
\end{equation}

Let $I\subset \{1,\dots, \dim(X)\}$ and let 
$$\mathrm{pr}_I:\bigoplus_{i=1}^{\dim(X)}L_i\ra \bigoplus_{i\in I}L_i$$ 
be the canonical projection. 

We claim that
there exists $J\subset \{1,\ldots, \dim(X)\}$ with $\sharp
J=\dim(Z)$ such that $i_J:=\mathrm{pr}_J\circ i$ is an
isomorphism, and hence $T_{Z_{\reg}}\simeq \bigoplus_{i\in
  J}L_i$.

 Indeed, there must be a subset $J\subset \{1,\ldots, \dim(X)\}$ 
 with 
 $\sharp J=\dim(Z)$ such that the map $i_J$ is injective, otherwise 
 $i$ in \eqref{seq} would not be injective. As $Z_{\reg}$ has no divisors, the degeneracy locus of $i_J$ is empty, and hence  $i_J$ is  an isomorphism, as claimed.

Thus we can see the map $i$ as a matrix 
$$A=(a_{ij}), i=1,\ldots,\dim(Z), j=1,\ldots,\dim(X)$$ 
with each $a_{ij}\in {\Hom}\left({L_i}\restrict{ Z_{\reg}}, {L_j}\restrict{Z_{\reg}}\right).$

Let $\pi:\tilde X\ra X$ be the universal cover of $X$. As  $Z_{\reg}$ has no divisors,  any morphism of line bundles of $Z_{\reg}$ is either zero or multiplication by a non-zero constant. As a consequence, 
the entries of the matrix $A$ are all constant, and hence the image of the bundle morphism 
$$T_{\pi^{-1}(Z_{\reg})}\ra {T_{\tilde{X}}}_{\vert \pi^{-1}(Z_{\reg})}$$
 is the vector subspace generated by the vectors $\{f_i\}, i\in J$, given by
$$f_i=\sum_{j=1}^{\dim(X)}a_{ij}\frac{\partial}{\partial z_j}.$$
In particular, the preimage $\pi^{-1}(Z_{\reg})$ of $Z_{\reg}$ on the universal cover $\tilde{X}=\H^s\times \C^t$ of $X$ is locally an open subset of an affine subspace ${\mathcal A}$ of $\C^{s+t}$ where the direction of ${\mathcal A}$ is spanned by  $$f_i=\sum_{j=1}^{\dim(X)}a_{ij}e_j$$ where $\{e_j\}, j=1,\dots, \dim(X)$ is the canonical basis in $\C^{s+t}.$

It follows that any connected component   $\tilde{Z}_{\reg}\subset\pi^{-1}(Z_{\reg})$   is contained as an open subset of an affine subspace ${\mathcal A}_{\tilde{Z}_{\reg}}$ of the same dimension, hence any connected component of  $\tilde{Z}\subset\pi^{-1}(Z)$ is also contained into such an affine subspace. We derive that  $\pi^{-1}(Z)$ is in fact smooth, and hence $Z$ is smooth, too.

Moreover, since $\tilde{Z}$ is closed in $\H^s\times \C^t$, the following equality holds:
$$\tilde{Z}={\mathcal A}_{\tilde{Z}} \cap \left(\H^s\times \C^t\right).$$

\hfill

\remark \label{_subvarieties_flat_Remark_}
Notice that this observation also proves \ref{_flat_subva_Theorem_}.
\endproof

\hfill



Now fix  a conected component $\tilde{Z}\subset \pi^{-1}(Z)$.  Then  
$$Z=\tilde{Z}/\Stab(\tilde{Z})$$ where 
$$\Stab(\tilde{Z})=\{\gamma\in \pi_1(X)\, \vert\,  \gamma(\tilde{Z})=\tilde{Z}\}.$$
Analysing the structure of the group $\Stab(\tilde{Z})$ will eventually lead to a contradiction. In the first place, observe that  $\Stab(\tilde{Z})$ cannot consist of translations only, since otherwise $Z$ would be a torus\footnote{A bit more details are needed here..}, contradicting with \ref{tori}. 

\hfill

Fix then $\gamma\in \Stab(\tilde{Z})$ which is not a translation and let $R_u$ be the linear map induced by $\gamma$.  The direction $\vec{\tilde{Z}}$ of the affine subspace $\tilde{Z}$ is then left invariant by $R_u$. Since $R_u$ is diagonal,  either  $\vec{\tilde{Z}}$ has a basis among the $\{e_1,\ldots, e_n\}$ 
or (at least) two of the eigenvalues of $R_u$ are equal. The second case is excluded by the assumption on $M$,
so there exists a subset $I\subset \{1,\dots, n\}$ such that $\vec{\tilde{Z}}=\oplus_{i\in I}\C e_i.$
It follows that for any $j\in I':=\{1,\dots, n\}\setminus I$ there exists constants $c_j\in \C$ such that $\tilde{Z}$ is given by the equations $z_j=c_j,\forall j\in I'.$ 

Let now $\gamma\in \Stab(\tilde{Z})$ be arbitrary. Note that $\gamma$ cannot be a translation by some $(\sigma_1(a),\dots\sigma_n(a))$, since then for any $j\in I'$ we would have $c_j=c_j+\sigma_j(a)$, yielding $\sigma_j(a)=0$, and hence $a=0.$ This implies that 
any nontrivial $\gamma\in \Stab(\tilde{Z})$ is of the form
$$\gamma(z_1,\dots, z_n)=(\sigma_1(u)z_1+\sigma_1(a),\dots, \sigma_n(u)z_n+\sigma_n(a)),$$
for some $u\in U, u\not=1$ and for some $a\in A_M.$

But since $\gamma\in \Stab(\tilde{Z})$ we see that for any $j\in J$ we have 
$$c_j=\sigma_j(u)c_j+\sigma_j(a),$$
and thus
$$c_j=\frac{\sigma_j(a)}{1-\sigma_j(u)},\quad \text{for all}\,  j\in J.$$
The point  $P_0\in \C^s\times \C^t$,
$$P_0=\left(\frac{\sigma_1(a)}{1-\sigma_1(u)}, \dots, \frac{\sigma_n(a)}{1-\sigma_n(u)}\right),$$
is thus fixed  by all $\gamma\in \Stab(\tilde{Z}).$
But then, after changing the coordinates in $\tilde{Z}$ via 
$$z_i\mapsto z_i-\frac{\sigma_i(a)}{1-\sigma_(u)},\quad i=1,\ldots, n,$$ 
 we see that 
$\Stab(\tilde{Z})$ acts on $\tilde{Z}$ by linear diagonal transformations.
This means that $\tilde{Z}$ has a compact quotient under the action of a group of diagonal transformations. But this is easily seen to be impossible since on one hand, if a free abelian group $G$ of linear diagonal transformations acts discretely then its rank is $1$ (look at the orbit through $G$ of any point), while since $\tilde{Z}$ is a contractible manifold of real dimension at least $2$, the rank of any free abelian group acting cocompactly  on it must equal its real dimension (cf \cite{ce}, Application 3, pp 357 and \cite{br} Example 5, pp 185).

This contradiction completes the  proof of \ref{main}. \endproof

\hfill

\remark The condition that ``any element $u\in U, u\not=1$
is a primitive element for $K$'' is satisfied by a wide
class of choices for $K$ and $U.$ For instance, if $K$ is
a number field of prime degree over $\Q$ then {\em any}
choice of the admissible group of units will satisfy this
condition.

\hfill

\remark Although the  \ref{expl2}  may suggest that the condition  that any $u\in U, u\not=1$ is a primitive element may be equivalent to the fact that the OT-manifold has no proper subvarieties, this is not entirely correct. There are cases when this condition is not satisfied, but the OT-manifold still  has no proper complex subvarieties. For instance, if $K$ is a number field with a single complex place ($t=1$), the OT-manifold $X$ has no proper complex subvarities by \cite{ov_mrl_11}. But for such a number field $K$ we see that the rank of the group of units of $K$ is
$s+1-1=s$ so if the number field $K$ contains a proper subfield $\Q{\subset}L \subset K$ then for any choice of the admissible group of units $U$ there are elements in $U$ which belong to $L,$ hence not all elements in $U$ are primitive elements for $K.$

\hfill

{\bf Acknowledgment:} We thank  
\c Stefan Papadima for the references \cite{ce} and \cite{br} and
Cezar Joi\c{t}a for the
references concerning  the theorem of Remmert-Stein.
L.O. and M.V. thank
Higher School of Economics, Moscow, and University of
Bucharest for facilitating mutual visits while working at
this paper.

{\small

\noindent {\sc Liviu Ornea\\
University of Bucharest, Faculty of Mathematics, \\14
Academiei str., 70109 Bucharest, Romania}, and:\\
{\sc Institute of Mathematics "Simion Stoilow" of the Romanian
Academy,\\
21, Calea Grivitei Str.
010702-Bucharest, Romania\\
\tt lornea@fmi.unibuc.ro, \ \  Liviu.Ornea@imar.ro

\hfill

\noindent {\sc Misha Verbitsky\\
            {\sc Instituto Nacional de Matem\'atica Pura e
              Aplicada (IMPA) \\ Estrada Dona Castorina, 110\\
Jardim Bot\^anico, CEP 22460-320\\
Rio de Janeiro, RJ - Brasil}\\
also:\\
{\sc Laboratory of Algebraic Geometry,\\
National Research University Higher School of Economics,\\
Department of Mathematics, 6 Usacheva street, Moscow, Russia.}\\

\hfill

\noindent {\sc Victor Vuletescu\\
University of Bucharest, Faculty of Mathematics, \\14
Academiei str., 70109 Bucharest, Romania.}\\
\tt vuli@fmi.unibuc.ro
}}
\end{document}